\documentclass{article}
\usepackage{amsmath,amsthm,amssymb,url}

\newtheorem{theorem}{Theorem}
\newtheorem{lemma}{Lemma}

\allowdisplaybreaks

\begin{document}
\centerline{\textbf{\Large On the irreducibility of a truncated
binomial expansion}}
\vskip 20pt
\centerline{\large by Michael Filaseta, Angel Kumchev and Dmitrii
V. Pasechnik}

\renewcommand{\thefootnote}{}
\footnote{
\vskip 0pt
The first author was supported by the National Science Foundation
during research for this paper.

The third author
was supported by NWO grant 613.000214.
Parts of the work were completed
while the third author was supported by MSRI and by
DFG grant SCHN-503/2-1, while he held positions at MSRI, resp. at
CS Dept., Uni. Frankfurt.}

\vskip 30pt
\section{Introduction}
For positive integers $k$ and $n$ with
$k \le n-1$, define
\[
P_{n,k}(x) = \sum_{j=0}^{k} \binom{n}{j} x^{j}.
\]
In the case that $k = n-1$, the polynomial $P_{n,k}(x)$ takes the form
\[
P_{n,n-1}(x) = (x+1)^{n} - x^{n}.
\]
If $n$ is not a prime, $P_{n,n-1}(x)$ is reducible over $\mathbb Q$.
If $n = p$ is prime, the
polynomial $P_{n,n-1}(x) = P_{p,p-1}(x)$ is irreducible as
Eisenstein's criterion applies
to the reciprocal polynomial $x^{p-1} P_{p,p-1}(1/x)$. This note
concerns the
irreducibility of $P_{n,k}(x)$ in the case where $1 \le k \le n-2$.
Computations for $n \le 100$ suggest that in this case $P_{n,k}(x)$
is always irreducible.
We will not be able to establish this but instead give some results
which give further
evidence that these polynomials are irreducible.

The problem arose during the 2004 MSRI program on ``\textit{Topological
aspects
of real algebraic geometry}'' in the context of work of Inna Scherbak
in her investigations of the Schubert calculus in Grassmannians. She
had observed
that the roots of any given $P_{n,k}(x)$ are simple. This follows
from the identity
\[
P_{n,k}(x) - (x+1) \dfrac{P'_{n,k}(x)}{n} = \binom{n-1}{k} x^k.
\]
She then asked whether, for a fixed positive integer $n$, the various
$n(n-1)/2$ roots of $P_{n,k}(x)$ for $1 \le k \le n-1$ are distinct.
We will not resolve this problem but our methods imply
that for each positive integer $n$, almost all of the roots are
distinct. In other words, the number of distinct roots is $\sim n^{2}/2$ as
$n$ tends to infinity. We note that since the initial writing of
this paper, Inna Scherbak \cite{is}
has resolved the original question concerning the roots of these
polynomials as a
consequence of her work in representation theory.

Before closing this introduction, we mention that these same
polynomials have recently
arisen in the context of work by Iossif V.~Ostrovskii \cite{o}. In
particular, he finds a solution
to a problem posed by Alexandre Eremenko on the distribution of the
zeroes of $P_{n,k}(x)$
as $k$ and $n$ tend to infinity with $k/n$ approaching a limit
$\alpha \in (0,1)$.

\section{The Results}
Our methods apply to a wider class of polynomials than the
$P_{n,k}(x)$'s alone, so we begin by
recasting the problem in a more general setting.
For $a$ and $b$ nonnegative integers with $a \le b$, the identity
\begin{equation}\label{eq1}
\sum_{j=0}^{a} \binom{b}{j} (-1)^{j}
= \binom{b-1}{a} (-1)^{a}
\end{equation}
is easily established by induction on $a$.
We deduce that
\begin{align*}
P_{n,k}(x-1) &= \sum_{j=0}^{k} \binom{n}{j} (x-1)^{j}
= \sum_{j=0}^{k} \binom{n}{j} \sum_{i=0}^{j} \binom{j}{i} (-1)^{j-i}
x^{i} \\[5 pt]
&= \sum_{i=0}^{k} \sum_{j=i}^{k} \binom{n}{j} \binom{j}{i} (-1)^{j-i} x^{i}
= \sum_{i=0}^{k} \sum_{j=i}^{k} \binom{n}{i} \binom{n-i}{j-i}
(-1)^{j-i} x^{i} \\[5 pt]
& = \sum_{i=0}^{k} \binom{n}{i} \sum_{j=0}^{k-i} \binom{n-i}{j} (-1)^{j}
x^{i}
= \sum_{i=0}^{k} \binom{n}{i} \binom{n-i-1}{k-i} (-1)^{k-i} x^{i},
\end{align*}
\pagestyle{plain}
where the last equality makes use of (\ref{eq1}).
For $0 \le j \le k$, we define
\begin{equation*}
c_{j} = \binom{n}{j} \binom{n-j-1}{k-j} (-1)^{k-j}
= \dfrac{(-1)^{k-j} n (n-1) \cdots (n-j+1) (n-j-1) \cdots (n-k+1)
(n-k)}{j! (k-j)!}
\end{equation*}
so that $P_{n,k}(x-1) = \sum_{j=0}^{k} c_{j} x^{j}$.
We are interested in the irreducibility of $P_{n,k}(x)$. A necessary
and sufficient
condition for $P_{n,k}(x)$ to be irreducible is for $P_{n,k}(x-1)$ to
be irreducible,
so we restrict our attention to establishing irreducibility results
for the polynomials
$\sum_{j=0}^{k} c_{j} x^{j}$.

For our results, we consider
\begin{equation}\label{neweq2}
F_{n,k}(x) = \sum_{j=0}^{k} a_{j} c_{j} x^{j},
\end{equation}
where $a_{0}, a_{1}, \dots, a_{k}$ denote integers, each having all
of its prime factors $\le k$.
In particular, none of the $a_{j}$ are zero.
Observe that if $F_{n,k}(x)$ is irreducible for all such $a_{j}$,
then necessarily $P_{n,k}(x)$
is irreducible simply by choosing each $a_{j} = 1$. Another
interesting choice for $a_{j}$ is
$a_{j} = (-1)^{k-j} j! (k-j)!$. As $F_{n,k}(x)$ is irreducible if and
only if
$((n-k-1)!/n!) \cdot x^{k} F_{n,k}(1/x)$ is irreducible, the
irreducibility of $F_{n,k}(x)$ will imply
the irreducibility of
\[
\dfrac{1}{n-k} + \dfrac{x}{n-k+1} + \dfrac{x^{2}}{n-k+2} + \cdots +
\dfrac{x^{k}}{n}.
\]
In particular, if $n = k+1$, these polynomials take a nice form. It
is possible to show, still
with $n = k+1$, that these polynomials are irreducible for every
positive integer $k$. The
idea is to use Newton polygons with respect to two distinct primes in
the interval
$((k+1)/2,k+1]$. For $k \ge 10$, it is known that such primes exist
(cf.~\cite{r}). As this is not the
focus of the current paper, we omit the details.

Let $N$ be a positive integer.
The number of integral pairs $(n,k)$ with $1 \le n \le N$ and $1 \le
k \le n-2$ is
\[
\sum_{n \le N} (n-2) \sim \dfrac{N^{2}}{2}.
\]
Our first result is that the number of possible reducible polynomials
$F_{n,k}(x)$
with $n \le N$ and $1 \le k \le n-2$ is small by comparison. More
precisely, we show
the following.

\begin{theorem}\label{thm1}
Let $\varepsilon > 0$, and let $N$ be a positive integer. For each
integral pair $(n,k)$ with
$1 \le n \le N$
and $1 \le k \le n-2$, consider the set $S(n,k)$ of all polynomials
of the form (\ref{neweq2})
where $a_{0}, a_{1}, \dots, a_{k}$ denote arbitrary integers, each
having all of its prime factors $\le k$.
The number of such pairs $(n,k)$ for which there exists a
polynomial $f(x) \in S(n,k)$ that is reducible is
$O(N^{23/18 + \varepsilon})$.
\end{theorem}

Under the assumption of the Lindel\"of hypothesis, a result of Gang
Yu \cite{gy}
can be used to improve our estimate for the number of exceptional pairs
$(n, k)$ to $O(N^{1 + \varepsilon})$.
A further improvement to $O(N \log^{3} N)$ is possible under the Riemann
hypothesis by a classical result of Atle Selberg \cite{as}.

Based on the main result in \cite{bhp}, one can easily modify
our approach to show that for each positive integer
$n$, there are at
most $O(n^{0.525})$ different positive integers $k \le n-2$ for which
$S(n,k)$ contains a
reducible polynomial. This implies the remark in the introduction
that for fixed $n$,
the number of distinct roots of $P_{n,k}(x)$ for $k \le n-2$ is $\sim
n^{2}/2$ as
$n$ tends to infinity.

Our second result is an explicit criterion for the irreducibility of
$F_{n,k}(x)$.

\begin{theorem}\label{thm2}
If there is a prime $p > k$ that exactly divides $n(n-k)$, then
$F_{n,k}(x)$ is irreducible
for every choice of integers $a_{0}, a_{1}, \dots, a_{k}$ with each
having all of its prime factors $\le k$.
\end{theorem}

Theorem \ref{thm2} has a simple proof based on Eisenstein's criterion.
It implies, in particular, that if $n$ is a prime, then $P_{n,k}(x)$
is irreducible
for every $k \in \{ 1, 2, \dots, n-1 \}$. This then resolves the
problem of Scherbak
in the case that $n$ is a prime.

Our third and final result concerning the irreducibility of
$F_{n,k}(x)$ is as follows.

\begin{theorem}\label{thm3}
Let $k$ be a fixed integer $\ge 3$. There is an $n_{0} = n_{0}(k)$ such
that
if $n \ge n_{0}$, then $F_{n,k}(x)$ is irreducible for every choice of
integers
$a_{0}, a_{1}, \dots, a_{k}$ with each having all of its prime factors
$\le k$.
\end{theorem}

The value of $n_{0}(k)$ in this last result, being based on the
solutions to certain Thue equations,
can be effectively determined. The result is of added interest as
the proof of Theorem \ref{thm1}
relies on considering $k$ large. Thus, the proof of Theorem
\ref{thm1} gives no
information about the situation in Theorem \ref{thm3}, where $k$ is
fixed and $n$ is large.
In the case that $k=1$, the polynomials $F_{n,k}(x)$ are linear and,
hence, irreducible.
In the case that $k=2$, our approach does not apply; but we note that
the polynomials
$P_{n,2}(x)$ are easily seen to be irreducible for $n \ge 3$ as
$P_{n,2}(x)$ has imaginary
roots for such $n$.

\section{The Proofs}
\begin{proof}[Proof of Theorem \ref{thm1}]
Let $f(x) = \sum_{j=0}^{k} d_{j} x^{j} \in \mathbb Z[x]$ with $d_{k}
d_{0} \ne 0$.
In the argument, we will make use of the Newton polygon of $f(x)$
with respect to
a prime $p$. The Newton polygon of $f(x)$ with respect to $p$ can be
defined as the lower part of
the convex hull of the points $(j,\nu_{p}(d_{j}))$ where $0 \le j \le
k$ and $\nu_{p}(m)$ is defined
to be the integer $r$ satisfying $p^{r}|m$ and $p^{r+1}\nmid m$.
Thus, the Newton polygon
has its left most endpoint being $(0,\nu_{p}(d_{0}))$ and its right
most endpoint being
$(k,\nu_{p}(d_{k}))$. A theorem of Gustave Dumas \cite{d} asserts that,
for a fixed prime, the Newton polygon
of a product of two polynomials can be obtained by translating the
edges of the Newton polygons
of each of the polynomials. The endpoints on the translation of an
edge always occur at lattice
points. For the proof of Theorem \ref{thm1}, we will use a
specific consequence of this result:
If the lattice points along the edges of the Newton polygon of $f(x)$
with respect to $p$ consist of $(0,\nu_{p}(d_{0}))$,
$(k,\nu_{p}(d_{k}))$ and
only one additional lattice point, say at $(u, v)$, then either
$f(x)$ is irreducible or it is the product of an irreducible
polynomial of degree $u$ times an irreducible polynomial of degree
$k-u$.
We note that the lattice point
$(u, v)$ in this context need not be one of the points
$(j,\nu_{p}(d_{j}))$ (for example, consider
$f(x) = x^{2} + 4x + 4$ and $p = 2$).

Consider $n$ sufficiently large.  Let $p_{j}$ denote the $j$th prime, and
let $t$ be maximal such that $p_{t} < n$.  Denote by $\delta(n)$ the
distance from $n$ to $p_{t-1}$ so that $\delta(n) = n - p_{t-1}$.  Suppose
that $k$ satisfies $2 \delta(n) < k < n- \delta(n)$.  We show that in this
case, the polynomial $f(x) = F_{n,k}(x)$ is irreducible over $\mathbb Q$
(independent of the choices of $a_{j}$ as in the theorem).  First, we
explain why this implies our result.

For the moment, suppose that we have shown that
$F_{n,k}(x)$ is irreducible over $\mathbb Q$ for $n$ sufficiently large
and $2 \delta(n) < k < n- \delta(n)$.
Let $\rho(n) = p_{t-1}$ where $t$ is defined as above.
It follows that the number of pairs $(n,k)$ as in the theorem
for which there exists a reducible polynomial $f(x) \in S(n,k)$ is
\[
\ll \sum_{n \le N} \delta(n)
\ll \sum_{n \le N} (n - \rho(n))
\ll \sum_{2 < p_{t} < N} \sum_{\substack{n \le N\\ \rho(n) =
p_{t-1}}} (n - p_{t-1}).
\]
This last double sum can be handled rather easily by extending the range
on $n$ slightly (to the least prime that is $\ge N$).  It does not exceed
\[
\sum_{2 < p_{t} < N} \sum_{p_{t} < n \le p_{t+1}} (n - p_{t-1})
\le \sum_{2 < p_{t} < N}  \sum_{p_{t} < n \le p_{t+1}} (p_{t+1} - p_{t-1})
\le \sum_{2 < p_{t} < N} (  p_{t+1} - p_{t-1} ) ( p_{t+1}-p_{t} ).
\]
Setting $d_{j} = p_{j+1}-p_{j}$, we deduce from the arithmetic-geometric
mean
inequality that
\[
  (p_{t+1} - p_{t})(p_{t+1}  - p_{t-1})
= d_{t} (d_{t} + d_{t-1}) = d_{t}^{2} + d_{t} d_{t-1}
\le {\textstyle \frac 32 d_{t}^{2} + \frac 12 d_{t-1}^{2}}.
\]
Hence, the number of pairs $(n,k)$ as in the theorem
for which there exists a reducible polynomial $f(x) \in S(n,k)$ is
$\ll \sum_{p_{t} < N} d_{t}^{2}$.
A theorem of Roger Heath-Brown \cite{rhb} asserts that
\[
\sum_{p_{t} \le N} d_{t}^{2} \ll N^{23/18 + \varepsilon}.
\]
Therefore, Theorem \ref{thm3} follows provided we establish that
$F_{n,k}(x)$ is irreducible over $\mathbb Q$ for $n$ sufficiently large
and $2 \delta(n) < k < n- \delta(n)$.

Consider $n$ sufficiently large and $k$ an integer in the interval
$(2 \delta(n), n- \delta(n))$.  Let $p = p_{t}$ and $q = p_{t-1}$.
Note that both $p$ and $q$ are greater than $k$. We set $u$
and $v$ to be the positive integers satisfying $p = n-u$ and $q = n-v$.
Then $1 \le u < v = \delta(n) < k/2$.
Observe that the numerator of $c_{j}$ is the product of the integers
from $n-k$ to $n$ inclusive
but with the factor $n-j$ missing.  Also, the denominator of $c_{j}$
is not divisible by any prime
$> k$ and, in particular, by $p$ or by $q$.

We look at the Newton polygon of $f(x)$ with respect to $p$ and the
Newton polygon of $f(x)$ with respect to $q$. Note that
$\nu_{p}(n-u) = 1$ and, for each $j$, we have $\nu_{p}(a_{j}) = 0$.
Therefore, the Newton polygon of $f(x)$ with respect to $p$ consists
of two line segments, one from $(0,1)$ to $(u,0)$ and one from
$(u,0)$ to $(k,1)$. The
theorem of Dumas implies that if $f(x)$ is reducible, then it must be
an irreducible polynomial of degree $u$ times an irreducible
polynomial of degree $k-u$.
Similarly, by considering the Newton polygon of $f(x)$ with respect
to $q$, we deduce
that if $f(x)$ is reducible, then it is an irreducible polynomial of
degree $v$ times an
irreducible polynomial of degree $k-v$. Since $k - v >  \delta(n) > u$ and
$v \ne u$, we deduce that $f(x)$ cannot be reducible. Thus, $f(x)$
is irreducible.
Theorem \ref{thm1} follows.
\end{proof}

\begin{proof}[Proof of Theorem \ref{thm2}]
Eisenstein's criterion applies
to $x^{k} F_{n,k}(1/x)$ whenever there is a prime $p > k$ that
exactly divides $n$
(i.e., $p|n$ and $p^{2}\nmid n$). Hence, $F_{n,k}(x)$
is irreducible whenever such a prime exists.
Also, $F_{n,k}(x)$ itself satisfies
Eisenstein's criterion whenever there is a prime $p > k$ that exactly
divides $n-k$.
\end{proof}

\begin{proof}[Proof of Theorem \ref{thm3}]
As in the previous proofs, we work with
$f(x) = F_{n,k}(x)$ (where the $a_{j}$ are arbitrary integers
divisible only by primes $\le k$).
With $k$ fixed, we consider $n$ large
and look at the factorizations of $n$ and $n-k$.

\begin{lemma}\label{lemma1}
Let $p$ be a prime $> k$ and $e$ a positive integer for which
$\nu_{p}(n) = e$ or $\nu_{p}(n-k) = e$. Then each irreducible factor
of $f(x)$ has degree a multiple of $k/\gcd(k,e)$.
\end{lemma}

The proof of Lemma \ref{lemma1} follows directly by considering the
Newton polygon
of $f(x)$ with respect to $p$. It consists of one edge, and the
$x$-coordinates of the
lattice points along this edge
will occur at multiples of $k/\gcd(k,e)$. The theorem of Dumas
implies that the
irreducible factors of $f(x)$ must have degrees that are multiples of
$k/\gcd(k,e)$.

\begin{lemma}\label{lemma2}
Let $n'$ be the largest divisor of $n(n-k)$ that is relatively prime
to $k!$. Write
\[
n' = p_{1}^{e_{1}} p_{2}^{e_{2}} \cdots p_{r}^{e_{r}},
\]
where the $p_{j}$ denote distinct primes and the $e_{j}$ are positive
integers.
Let
\begin{equation}\label{eq3}
d = \gcd(k,e_{1}, e_{2}, \dots, e_{r}).
\end{equation}
Then the degree of each irreducible factor of $f(x)$ is a multiple of
$k/d$.
\end{lemma}

Note that in the statement of Lemma \ref{lemma2}, if $n' = 1$, then $d =
k$.
The proof of Lemma \ref{lemma2} makes use of Lemma \ref{lemma1}. Suppose
$m$ is the degree of an irreducible factor of $f(x)$. Then for each
$j \in \{ 1, 2, \dots, r \}$,
Lemma \ref{lemma1} implies there is an integer $b_{j}$ such that $m
e_{j} = k b_{j}$.
There are integers $x_{j}$ for which
\[
k x_{0} + e_{1} x_{1} + e_{2} x_{2} + \cdots + e_{r} x_{r} = d.
\]
Hence,
\[
m (d - k x_{0}) = m \big( e_{1} x_{1} + e_{2} x_{2} + \cdots + e_{r}
x_{r} \big)
= k \big( b_{1} x_{1} + b_{2} x_{2} + \cdots + b_{r} x_{r} \big).
\]
It follows that $m d$ is a multiple of $k$ so that $m$ is a multiple
of $k/d$ as claimed.

For the proof of Theorem \ref{thm3}, we define $d$ as in Lemma
\ref{lemma2} and consider three cases:
(i) $d = 1$, (ii) $d = 2$, and (iii) $d \ge 3$. In Case (i), Lemma
\ref{lemma2} implies $f(x)$ is irreducible. Case (ii) is more
difficult
and we return to it shortly. In Case (iii), there exist positive
integers $a$, $b$, $m_{1}$, and $m_{2}$
satisfying
\begin{equation}\label{eq4}
n = a m_{1}^{d}, \quad n-k = b m_{2}^{d}, \quad \text{ and } \quad a
\text{ and } b \text{ divide } \prod_{p \le k} p^{d-1}.
\end{equation}
As $d|k$ and $k$ is fixed, there are finitely many choices for $d$,
$a$ and $b$ as in (\ref{eq4}).
For each such $d$, $a$ and $b$, the possible values
of $n$ correspond to $a x^{d}$ given by solutions to the Diophantine
equation
\[
a x^{d} - b y^{d} = k.
\]
The above is a Thue equation, and it is well known that, since $d \ge
3$, it has finitely many solutions in integers $x$ and $y$
(see \cite{st}). It follows that there are finitely many integers
$n$ for which Case (iii) holds. Hence, for $n$ sufficiently
large, Case (iii) cannot occur. We are left with considering Case (ii).

For Case (ii), the definition of $d$ implies $k$ is even. As we
already have $k \ge 3$,
we deduce $k \ge 4$. Lemma \ref{lemma2} implies that
if $f(x)$ is reducible, then it factors as a product of two
irreducible polynomials each of
degree $k/2$. To finish the analysis for Case (ii), we make use of
the following.

\begin{lemma}\label{lemma3}
Let $f(x)$ be as above with $d$, as defined in (\ref{eq3}), equal to
$2$.ively
prime to $k!$.
Suppose $\nu_{p}(n'') = e$ where $p$ is a prime $> k$ and $e$ is a
positive integer.
If $f(x)$ is reducible, then $(k-1)|e$.
\end{lemma}

For the proof of Lemma \ref{lemma3}, we again appeal to the theorem of
Dumas.
Suppose first that $p|(n-1)$. Since $p > k$, we deduce that
$\nu_{p}(n-1) = e$.
The Newton polygon of $f(x)$ with respect to $p$ consists of two line
segments,
one from $(0,e)$ to $(1,0)$ and one from $(1,0)$ to $(k,e)$. Let $d'
= \gcd(k-1,e)$.  As $d = 2$, we deduce as above that $f(x)$ is a product
of two irreducible polynomials of degree $k/2$.  The fact that $f(x)$ has
just two irreducible factors implies by the theorem of Dumas
that one of these factors has degree
that is a multiple of $(k-1)/d'$ and the other has degree that is one more
than a multiple of $(k-1)/d'$.   We deduce that there are integers
$m$ and $m'$ such that
\[
  \dfrac{k-1}{d'} m = \dfrac{k}{2} \quad \text{ and } \quad
  \dfrac{k-1}{d'} m' + 1 = \dfrac{k}{2}.
\]
It follows that $(k - 1)/d'$ divides $1$, whence
\[
k-1 = d' = \gcd(k-1,e).
\]
We deduce that $(k-1)|e$.  A similar argument works in the case that
$p|(n-k+1)$.

To finish the analysis for Case (ii), we use Lemma \ref{lemma3} to
deduce that there are
positive integers $a'$, $b'$, $m_{3}$, and $m_{4}$ such that
\begin{equation}\label{eq5}
n-1 = a' m_{3}^{k-1}, \quad n-k+1 = b' m_{4}^{k-1}, \quad \text{ and
} \quad a' \text{ and } b' \text{ divide } \prod_{p \le k} p^{k-2}.
\end{equation}
As $k$ is fixed, there are finitely many choices for $a'$ and $b'$ as
in (\ref{eq5}).
For each of these, the possible values
of $n$ correspond to $a' x^{k-1} + 1$ determined by solving the
Diophantine equation
\[
a' x^{k-1} - b' y^{k-1} = k-2.
\]
As $k \ge 4$, the above is a Thue equation and has finitely many
solutions in integers $x$ and $y$.
Thus, there are finitely many integers $n$ for which Case (ii) holds.
Hence, for $n$ sufficiently
large, we deduce that $F_{n,k}(x)$ is irreducible, completing the
proof of Theorem \ref{thm3}.
\end{proof}

\vskip 15pt \footnotesize
\centerline{\hrulefill}
\vskip 5pt
\centerline{\hskip 40pt
\parbox[t]{6cm}{
Mathematics Department
\vskip 0pt \noindent
University of South Carolina
\vskip 0pt \noindent
Columbia, SC  29208
\vskip 0pt \noindent
USA
\vskip 5pt \noindent
Email: \url{filaseta@math.sc.edu}
\vskip 5pt \noindent
\url{http://www.math.sc.edu/~filaseta/}
}
\hfill
\parbox[t]{6cm}{
Department of Mathematics
\vskip 0pt \noindent
The University of Texas at Austin
\vskip 0pt \noindent
1 University Station, C1200
\vskip 0pt \noindent
Austin, TX 78712-0257
\vskip 0pt \noindent
USA
\vskip 5pt \noindent
Email: \url{kumchev@math.utexas.edu}
\vskip 5pt \noindent
\url{http://www.ma.utexas.edu/~kumchev/}
}
\hfill
\parbox[t]{6cm}{
Dept. E \& OR
\vskip 0pt \noindent
Tilburg University, P.O.~Box 90153
\vskip 0pt \noindent
5000 LE Tilburg
\vskip 0pt \noindent
The Netherlands
\vskip 5pt \noindent
Email: \url{d.v.pasechnik@uvt.nl}
}
}

\end{document}